\documentclass[10pt, reqno, a4paper]{amsart}
\usepackage{amsmath,amsthm,amssymb,amsfonts}

\usepackage{enumitem}
\usepackage{a4wide}
\usepackage{bm}
\usepackage{mathtools}
\usepackage{thmtools}
\usepackage[dvipsnames]{xcolor}
\usepackage[colorlinks=true]{hyperref}
\usepackage{graphicx}
\usepackage{orcidlink}

\theoremstyle{plain}
\declaretheorem[name=Proposition, numberwithin=section]{proposition}
\declaretheorem[name=Theorem, sibling=proposition]{theorem}
\declaretheorem[name=Lemma, sibling=proposition]{lemma}

\theoremstyle{definition}
\declaretheorem[name=Remark, sibling=proposition]{remark}

\DeclareMathOperator{\E}{\mathsf{E}}

\renewcommand{\P}{\mathsf{P}}

\DeclareMathOperator{\N}{\mathsf{N}}
\DeclareMathOperator{\R}{\mathsf{R}}

\DeclareMathOperator{\D}{\mathsf{D}}
\renewcommand{\L}{\mathsf{L}}
\renewcommand{\d}{\mathrm{d}}
\DeclareMathOperator{\Dom}{\mathrm{Dom}\,}

\DeclareMathOperator*{\plim}{\P-\lim\,}

\newcommand{\supnorm}[2]{\left| #1 \right|_{\infty, #2}}
\newcommand{\numberthis}{\addtocounter{equation}{1}\tag{\theequation}}

\begin{document}
\numberwithin{equation}{section}

\title[$1/H$-variation of divergence integral with respect to a Hermite process]{On the $1/H$-variation of the divergence integral with respect to a Hermite process}

\author{Petr \v Coupek\orcidlink{0000-0002-5360-6095}}
\address{Charles University, 
		Faculty of Mathematics and Physics, 
		Sokolovsk\' a 83, 
		Prague 8, 186 75, 
		Czech Republic. }
\email{coupek@karlin.mff.cuni.cz}

\author{Pavel K\v r\'i\v z\orcidlink{0000-0001-5773-9541}}
\address{Charles University, 
		Faculty of Mathematics and Physics, 
		Sokolovsk\' a 83, 
		Prague 8, 186 75, 
		Czech Republic.}
\email{kriz@karlin.mff.cuni.cz}

\author{Mat\v{e}j Svoboda}
\address{Charles University, 
		Faculty of Mathematics and Physics, 
		Sokolovsk\' a 83, 
		Prague 8, 186 75, 
		Czech Republic.}
\email{matej.svoby@gmail.com}

\begin{abstract}
In this paper, a divergence-type integral of a random integrand with respect to the Hermite process of order $k\in\N$ with Hurst parameter $H\in (1/2,1)$ is defined and it is shown that the integral is of finite $1/H$-variation.
\end{abstract}

\keywords{Hermite process, $p$-variation, divergence integral, Malliavin calculus.}

\subjclass[2020]{60G22; 60H07}


\maketitle


\section{Introduction}
The concept of the $p$-variation can be thought of as a natural descriptor of the roughness of various stochastic processes; and, in fact, many statistical procedures designed for fitting stochastic models of complex systems with randomness are based on the quadratic, or, more generally, $p$-variation; see, e.g., \cite{ChrTudVie09, ChrTudVie11, TudVie07} and many others. 

For $p\in (0,\infty)$, $T\in (0,\infty)$, $n\in\N$, the dyadic partition $\{t_i^n=Ti2^{-n}\}_{i=0}^{2^n}$ of interval $[0,T]$, and for a stochastic process $X=(X_t)_{t\in [0,T]}$, let 
	\[ V_{n,T}^p(X) := \sum_{i=0}^{2^n-1} |X_{t_{i+1}^n}-X_{t_i^n}|^p\]
and define the \emph{$p$-variation of process $X$ on the interval $[0,T]$} by
	\[ V_T^p(X) := \plim_{n\to\infty} V_{n,T}^p(X)\]
provided that the limit exists.

One of the basic stochastic processes for which the $p$-variation is known and commonly used is, besides the Wiener process, the fractional Brownian motion (abbr.\ FBM); see, e.g., \cite{DecUst99,ManNes68}. The FBM $B^H=(B^H_t)_{t\in [0,T]}$ with Hurst parameter $H\in (0,1)$ is a centered Gaussian process whose covariance function is given by
	\[ R_H(s,t) = \frac{1}{2} \left(s^{2H} + t^{2H} - |s-t|^{2H}\right), \quad s,t\in [0,T].\]
As such, the FBM can be shown to have a continuous modification; it has stationary increments that can be negatively correlated (if $H<1/2$), independent (if $H=1/2$), or positively correlated (if $H>1/2$); and it also allows to model the source of randomness that is with or without memory (if $H\neq 1/2$ or $H=1/2$, respectively). It is not surprising that its $p$-variation is related to the value of the Hurst parameter and it was first shown in \cite{Rog97} that there is the convergence
	\[V_{n,T}^{p}(B^H) \quad\xrightarrow[n\to\infty]{\P}\quad \begin{cases}
		0, &\quad pH>1,\\		
		T C_{H,1}, &\quad pH=1,\\
		\infty, &\quad pH<1,
		\end{cases} \]
where $C_{H,1}:=\E|B_1^H|^{\frac{1}{H}}$. For models with additive fractional noise, such result would be sufficient; however, for models with a more complicated (multiplicative) noise structure, the stochastic integral with respect to the FBM has to be considered. However, except for the case $H=1/2$, the $H$-FBM is not a semimartingale so that one cannot use the standard (It\^o) stochastic calculus and other approaches have to be considered. One of the notable approaches builds upon the results of the Malliavin calculus and uses the so-called Skorokhod (or the divergence-type) integral. This operator is defined as the $L^2$-adjoint of the Malliavin derivative (see, e.g., \cite[Definition 1.3.1]{Nua06} for the precise definition) and extends the definition of the standard It\^o integral with respect to the Wiener process to nonanticipative integrands (see, e.g., \cite[Proposition 1.3.11]{Nua06} for the precise relationship). In particular, the stochastic integral $\int_0^T u_s\delta B^H_s$ of a suitable stochastic integrand $u$ with respect to the FBM $B^H$ is then defined as the composition of the Skorokhod integral with a certain deterministic fractional transfer operator (see, e.g., \cite{AloMazNua01, CarCouMon03, DecUst99}, or \cite{BiaHuOksZha08} and the many references therein). With such definition, the following results on the $1/H$-variation of the stochastic integral is obtained in \cite[Theorem 4.4]{GueNua05} and \cite[Theorem 4.1]{EssNua15} for $H>1/2$ and $H<1/2$, respectively:
	\[ V_{n,T}^\frac{1}{H}\left(\int_0^\bullet u_s\delta B_s^H\right) \quad\xrightarrow[n\to\infty]{ L^1(\Omega) }\quad C_{H,1}\int_0^T |u_s|^\frac{1}{H}\d{s}. \]
	
The FBM is the only Gaussian self-similar process with stationary increments; see, e.g., \cite{EmbMae02}. There are, however, many other non-Gaussian self-similar processes with stationary increments \cite{EmbMae02, MaeTud12} that can be considered as alternatives to the FBM in situations where a Gaussian process seems to be an unsuitable model (as in, e.g.,  modeling industrial control error \cite{Dom15}). Among all the viable models, the family of Hermite process (abbr. HP) might provide an interesting choice. Any HP is a self-similar process with stationary increments (and as such, its covariance function is also $R_H$) that arises as a limit of suitably normalized sums of strongly correlated random variables in the non-central limit theorem \cite{DobMaj79, Taq79}. The closest non-Gaussian process to the FBM within this class is the so-called Rosenblatt process (abbr.\ RP) and as such, it has been given considerable attention in recent years. In particular, a thorough exposition of the construction and basic properties of the RP can be found in the survey article \cite{Taq11}, its in-depth analysis in \cite{Tud08}, and some of its finer properties in, e.g.,  \cite{AbrPip06,Alb98,CouOnd24,DawLoo22,GarTorTud12,KerNouSakVii21,Pip04}. Additionally, the behavior of solutions to stochastic (partial) differential equations (abbr.\ S(P)DE) with Rosenblatt noise has also been studied; see, e.g., \cite{BonTud11,CouMasOnd18,CouMasOnd22,SlaTud19,SlaTud19b,SlaTud19c}. Of course, the general class of HPs has also been considered in the literature and we refer to, e.g., \cite{Loo23,PipTaq17,Tud13,Tud23} for an excellent overview that contains many references, to \cite{Aya20,AyaHamLoo25,MorOod86} for some finer results, to \cite{CouHenSla24, CouHenSla25, MaeTud07, AssDieTud22, SlaTud19, SlaTud19c} for results on SDEs driven by the HP, and to \cite{AssTud20, ChrTudVie11, DouEsSTud20, LooTud25} for statistical inference for the HP and HP-driven SDEs.

As far as variations of the HP $Z^{H,k}=(Z^{H,k}_t)_{t\in [0,T]}$ of order $k\in\N$ with Hurst parameter $H\in (1/2,1)$ are concerned, there are some results on the (normalized) quadratic variation, see \cite{Tud13}, and it is also known that there is the convergence
	\begin{equation}
	\label{eq:pvar_RP} 
		V_{n,T}^{\frac{1}{H}}(Z^{H,k}) \quad \xrightarrow[n\to\infty]{\P}\quad C_{H,k}T
	\end{equation}
where $C_{H,k} := \E|Z_1^{H,k}|^{\frac{1}{H}}$; see \cite[Proposition 6]{AssTud20}. It is therefore natural to ask whether also the convergence 
	\[
	V_{n,T}^\frac{1}{H}\left(\int_0^\bullet u_s\delta Z_s^{H,k}\right) 
	\quad\xrightarrow[n\to\infty]{ L^1(\Omega) }\quad 
	C_{H,k}\int_0^T |u_s|^\frac{1}{H}\d{s} 
	\]
for a divergence-type integral with respect to the HP holds for a sufficiently rich class of stochastic integrands $u=(u_t)_{t\in [0,T]}$. In the article, we first construct the integral and then give the affirmative answer to this question.

The construction of the divergence-type integral with respect to a general HP is analogous to the construction of the divergence-type integral with respect to the FBM and RP (see, e.g., \cite{AloMazNua01, BiaHuOksZha08, CarCouMon03, DecUst99} and \cite{CouDunDun22, Tud08}, respectively) but it appears not to have been given in the literature so far. The general proof strategy of \autoref{thm:main} is motivated by the proof of \cite[Theorem 4.4]{GueNua05} but the fact that we consider an HP of an arbitrary order requires additional non-trivial novel results on dense subspaces of the Sobolev--Watanabe spaces and a generalization of the pull-out rule for the multiple divergence integral that are of independent interest.

\emph{Organization of the article.} The article is organized as follows. In \autoref{sec:prelim}, we recall some basic notions from the Malliavin calculus and the definition of HPs. Subsequently, we continue with the construction of the divergence-type integral with respect to an HP. The main result is given in \autoref{sec:main} and in particular in \autoref{thm:main}. The auxiliary results on the Malliavin calculus and several technical lemmas needed for the proof of \autoref{thm:main} are given in \autoref{sec:app}.

\section{Preliminaries}
\label{sec:prelim}
We begin by recalling some concepts from the Malliavin calculus, the definition of the HP, and the definition of the divergence-type integral with respect to the HP. We refer the reader to the excellent monographs \cite{NouPec12} or \cite{Nua06} for a thorough exposition of the Malliavin calculus.

\subsection{Conventions and notation.} Throughout the article, we write $A\lesssim B$ if that there exists a
finite positive constant $c$ such that $A \leq cB$. The constant is independent of the variables which appear in the expressions $A$ and $B$ and it can change from line to line. This symbol is used when the exact value of the constant is not important for the exposition. Similarly, we also use the symbol $A\eqsim B$ if both $A\lesssim B$ and $B\lesssim A$ hold and the symbol $A\propto B$ if there exists a finite positive constant $c$ such that $A=cB$. 

For $n\in\N$, we write $C_b^\infty(\R^n)$ for the space of infinitely differentiable functions $f:\R^n\to\R$ such that $f$ and all its derivatives $\partial^k_{i_1,\ldots,i_k}f$ for $i_1,\ldots, i_k\in \{1,\ldots, n\}$ and $ k\in\N$, are bounded.

\subsection{Some notions from the Malliavin calculus}
Let $T\in (0,\infty)$ and let $(\Omega,\mathcal{F},\P)$ be a complete probability space with a standard Wiener process $W=(W_t)_{t\in [0,T]}$ defined on it. We assume that $\sigma$-algebra $\mathcal{F}$ is generated by process $W$. Let $I$ denote the \emph{Wiener--It\^o integral} (i.e.\ the unique extension of the linear isometry $\bm{1}_{(0,t]} \mapsto W_t$ from step functions on $[0,T]$ to $L^2(\Omega)$ to a linear isometry from $L^2([0,T])$ to $L^2(\Omega)$) and let $\mathcal{S}_b$ be the \emph{space of smooth bounded random variables} (i.e.\ the space of random variables $F=f(I(h_1),\ldots, I(h_n))$ for $n\in\N$, $f\in C^\infty_b(\R^n)$, and $h_1,\ldots,h_n\in L^2([0,T])$). For $k\in\N$ and $p\in [1,\infty)$, let $\D^{k,p}$ denote the natural \emph{domain of the $k$\textsuperscript{th} Malliavin derivative} (i.e.\ the closure of space $\mathcal{S}_b$ with respect to the norm 
	\[ \|F\|_{\D^{k,p}} := \left(\E |F|^p + \E\|D^1F\|_{L^2([0,T])}^p + \ldots + \E\|D^k F\|_{L^2([0,T]^k)}^p\right)^\frac{1}{p};\]
sometimes also called the Sobolev--Watanable space) and let $D^k$ denote the $k$\textsuperscript{th} \emph{Malliavin derivative} itself (i.e.\ the closed extension of the operator $D^k: L^p(\Omega)\supset \mathcal{S}_b\to L^p(\Omega;L^2([0,T]^k))$ defined for $F\in\mathcal{S}_b$ by 
	\[ D^kF := \sum_{i_1=1}^n\ldots \sum_{i_k=1}^n (\partial_{i_1, \ldots, i_k}^k f)(I(h_1),\ldots, I(h_n)) h_{i_1}\otimes\ldots \otimes h_{i_k},\]
to an operator $\D^{k,p}\to L^p(\Omega;L^2([0,T]^k)$). We also denote by $\Dom \delta^k$ the \emph{domain of the $k$\textsuperscript{th} order divergence operator} (i.e.\ the space of stochastic processes $u\in L^2(\Omega;L^2([0,T]^k))$ for which there exists $c\in (0,\infty)$ such that the inequality
	\[ |\langle D^kF,u\rangle_{L^2(\Omega;L^2([0,T]^k)}|\leq c \|F\|_{L^2(\Omega)}\]
holds for every $F\in\mathcal{S}_b$) and by $\delta^k$ the \emph{$k$\textsuperscript{th} order divergence operator} itself (i.e.\ the operator $\delta^k: \Dom \delta^k\to L^2(\Omega)$ defined for $u\in \Dom \delta^k$ as the unique element of $L^2(\Omega)$ for which
	\[ \langle F,\delta^k(u)\rangle_{L^2(\Omega)} = \langle u, D^kF\rangle_{L^2(\Omega; L^2([0,T]^k)}\]
holds for every $F\in\mathcal{S}_b$). For simplicity of notation we write $\D^{0,p}$ for $L^p(\Omega)$, $D^0$ for the identity operator on $L^p(\Omega)$, $\Dom\delta^0$ for $L^2(\Omega)$, and $\delta^0$ for the identity on $L^2(\Omega)$. We also identify the space $L^2([0,T]^0)$ with $\R$. Additionally, instead of $D^1$ and $\delta^1$, we write simply $D$ and $\delta$, respectively.

\subsection{Definiton of Hermite processes}
For our analysis, we consider the following finite-time interval representation of the HP; see \cite[Definition 7]{NouNuaTud10}, \cite[Theorem 1.1]{PipTaq10}, or \cite[Theorem 3.1]{Tud13}. For $H\in (1/2,1)$, $k\in\N$, and $t\in [0,T]$, define 
the kernel $L^{H,k}_t: [0,T]^k\to\R$ by
	\[ 
		L^{H,k}_t(x) 
			:= c_{H,k} \left[\int_{x_1\vee \ldots\vee x_k}^t \prod_{i = 1}^k \left(\frac{u}{x_i}\right)^{\frac{1}{2}-\frac{1-H}{k}} (u-x_i)^{-\left(\frac{1}{2}+\frac{1-H}{k}\right)}\d{u}\right] \bm{1}_{(0,t)^k}(x)
	\]
for $x=(x_1,\ldots,x_k)\in (0,T]^k$ where $c_{H,k}$ is defined by
	\begin{equation}
	\label{eq:cHk}
	c_{H,k} := \left(\frac{H(2H-1)}{k!\mathrm{B}\left(\frac{1}{2}-\frac{1-H}{k},\frac{2(1-H)}{k}\right)^k}\right)^\frac{1}{2}
	\end{equation}
in which $\mathrm{B}(\cdot\,,\cdot)$ is the Beta function. We have $L_t^{H,k}\in L^2([0,T]^k)$ for every $t\in [0,T]$ and we can therefore define the \emph{Hermite process} (of order $k$ with Hurst index $H$) $Z^{H,k}=(Z_t^{H,k})_{t\in [0,T]}$ by 
	\[
		Z_t^{H,k} := \delta^k(L^{H,k}_t), \quad t\in [0,T],
	\]
where $\delta^k$ is the $k$\textsuperscript{th} order divergence operator defined in the previous section but since the function $L^{H,k}_t$ is deterministic, this operator coincides with the Wiener--It\^o multiple integral of order $k$.
	
\subsection{Stochastic integral with respect to a Hermite process}
It immediately follows from \eqref{eq:pvar_RP} that no HP is a semimartingale; and as a consequence, a stochastic integral with respect to a HP cannot be understood in the classical It\^o sense. One approach to define the stochastic integral is via the Malliavin calculus.

In particular, let $E$ be a normed linear space and denote by $\mathcal{E}([0,T],E)$ the vector space of $E$-valued step functions defined on interval $[0,T]$; i.e.\ each $g\in\mathcal{E}([0,T];E)$ can be written as 
	\begin{equation}
	\label{eq:elementary_g} 
		g = \sum_{i=0}^{n-1} G_i\bm{1}_{[t_i,t_{i+1})}
	\end{equation}
for some $n\in\N$, $\{G_i\}_{i=0}^{n-1}\subset E$, and a partition $\{0=t_0\leq \ldots \leq t_n=T\}$ of the interval $[0,T]$. For $H\in (1/2,1)$ and $k\in\N$, define also
	\begin{equation}
	\label{eq:LHk}
		(L^{H,k}g)(x) := c_{H,k}\int_{0}^T g(u) \prod_{i = 1}^k \left(\frac{u}{x_i}\right)^{\frac{1}{2}-\frac{1-H}{k}} (u-x_i)_+^{-\left(\frac{1}{2}+\frac{1-H}{k}\right)}\d{u}
	\end{equation}
for $g\in\mathcal{E}([0,T],E)$ and $x=(x_1,\ldots, x_k)\in (0,T]^k$. In the above expression, $c_{H,k}$ is the constant defined by \eqref{eq:cHk}, $(u)_+:=\max\{0,u\}$, and the integral should be interpreted as the corresponding sum of elements in space $E$. There is the following boundedness result:
 
\begin{lemma}
\label{lem:bddness_of_LH}
For every $H\in (1/2,1)$, $k\in\N$, and $g\in \mathcal{E}([0,T],E)$, there is the estimate
	\[
		\|L^{H,k}g\|_{L^2([0,T]^k,E)} \lesssim \| g\|_{L^\frac{1}{H}([0,T],E)}.
	\]
\end{lemma}

\begin{proof}
By using the Fubini theorem, there is the estimate
	\[
		 \int_{[0,T]^k} \left\|(L^{H,k}g)(x)\right\|_E^2\,\d{x} \lesssim \int_0^T\int_0^T \|g(u)\|_E\|g(v)\|_E\, K(u,v)\,\d{u}\d{v}
	\]
where 
	\[
		K(u,v) := \left((uv)^{\frac{1}{2}-\frac{1-H}{k}}\int_0^{u\wedge v} x^{-2\left(\frac{1}{2}-\frac{1-H}{k}\right)} (u-x)^{-\left(\frac{1}{2}+\frac{1-H}{k}\right)}(v-x)^{-\left(\frac{1}{2}+\frac{1-H}{k}\right)} \d{x}\right)^k.
	\]
If $v< u$, it follows by substitution $y= \frac{u-x}{v-x}$ and then by substitution $z=\frac{u}{vy}$ that 
	\[ \int_0^v x^{-2\alpha} (u-x)^{\alpha-1}(v-x)^{\alpha-1}\d{x} = \mathrm{B}(\alpha,1-2\alpha)(uv)^{-\alpha} (u-v)^{2\alpha-1} \]
holds for any $\alpha\in (0,1/2)$ (and similarly in the case $v>u$). By using this fact, we obtain
	\[ 
		K(u,v) = \mathrm{B}\left(\frac{1}{2}-\frac{1-H}{k}, \frac{2(1-H)}{k}\right)^k |u-v|^{2H-2}, \quad u\neq v,
	\]
so that 
	\[
		 \| L^{H,k}g\|_{L^2([0,T]^k,E)}^2 \lesssim \int_0^T\int_0^T \|g(u)\|_E\|g(v)\|_E |u-v|^{2H-2}\d{u}\d{v}.
	\]
The claim now follows by standard arguments based on the symmetry of the integrand, the H\"older inequality with $p=1/H$, and the mapping properties of the one-sided Riemann--Liouville fractional integral from \cite[Theorem 3.5]{SamKilMar93} as in, e.g., \cite[Theorem 1.1]{MemMisVal01}. 
\end{proof}

Denote by $\mathcal{S}_{\infty,b}$ the space of random variables $F$ such that $F=f(I(h_1),\ldots, I(h_n))$ for some $n\in\N$, $f\in C_b^\infty(\R^n)$, and $h_1,\ldots, h_n\in C([0,T])$. The reason for introducing class $\mathcal{S}_{\infty,b}$ is that every $F$ of the above form satisfies 
	\[ 
		\sup_{\omega\in\Omega} \sum_{\ell=0}^k \sup_{x\in [0,T]^\ell} |(D^\ell F)(\omega,x)| <\infty
	\]
for every $k\in\N$ since $f\in C_b^\infty(\R^n)$ and since $h_1,\ldots h_n$ are all continuous, and therefore bounded, functions on $[0,T]$. This property will be crucial in the proof of the main result of the article given in \autoref{thm:main}. Let also $\mathcal{E}_{\infty,b}:=\mathcal{E}([0,T],\mathcal{S}_{\infty,b})$ be the space of $\mathcal{S}_{\infty,b}$-valued elementary processes defined on the interval $[0,T]$, i.e.\ the space of all processes $g$ of the form \eqref{eq:elementary_g} where $\{G_i\}_{i=0}^{n-1} \subset \mathcal{S}_{\infty,b}$.
For $H\in (1/2,1)$, $k\in\N$, and $g\in \mathcal{E}_{\infty,b}$, we define the \emph{Skorokhod (or the divergence-type) integral of $g$ with respect to the HP $Z^{H,k}$} by 
	\[ 
		\int_0^T g_s\delta Z_s^{H,k} := \delta^k (L^{H,k} g).
	\]

\begin{remark}
Note that if $H\in(1/2,1)$, $k\in\N$ and if $g\in\mathcal{E}_{\infty,b}$ is of the form \eqref{eq:elementary_g} with $\{G_i\}_{i=0}^{n-1}$ deterministic, we obtain the equality
	\[ \int_0^Tg_s\delta Z_s^{H,k} = \sum_{i=0}^{m-1}G_i(Z^{H,k}_{t_{i+1}}-Z^{H,k}_{t_i})\]
so that the above-defined divergence-type integral extends the Wiener-- It\^o integral with respect to the HP $Z^{H,k}$ to stochastic integrands. 
\end{remark}

We have the following boundedness result that is a particular case of the claim in \autoref{lem:boundedness_of_int_appendix}. For its purposes, let us denote $\L^{k,\frac{1}{H}}:= L^\frac{1}{H}([0,T],\D^{k,\frac{1}{H}})$ for $H\in (1/2,1)$ and $k\in\N$. 

\begin{lemma}
\label{lem:bddness_of_integral}
For every $H\in (1/2,1)$, $k\in\N$, and $g\in \mathcal{E}_{\infty,b}$, there is the estimate
	\begin{equation}
	\label{eq:bound_for_int}
	\left\|\int_0^T g_s\delta Z_s^{H,k}\right\|_{L^\frac{1}{H}(\Omega)}\lesssim \|g\|_{\L^{k,\frac{1}{H}}}.
	\end{equation}
\end{lemma}

There is also the following density result that is a particular case of the claim in \autoref{lem:density_of_E_appendix}.

\begin{lemma}
\label{lem:density_of_E}
The space $\mathcal{E}_{\infty,b}$ is dense in the space $\L^{k,\frac{1}{H}}$ for every $H\in (1/2,1)$ and $k\in\N$.
\end{lemma}

It follows by the above two lemmas that for every $H\in (1/2,1)$ and $k\in\N$, the operator $\int_0^T (\,\bullet\,)\delta Z_s^{H.k}$ can be extended to a bounded linear operator from $\L^{k,\frac{1}{H}}$ to $L^\frac{1}{H}(\Omega)$ for which \eqref{eq:bound_for_int} holds. We denote the extension by the same symbol and for $g\in \L^{k,\frac{1}{H}}$, we call the random variable $\int_0^T g_s\delta Z_s^{H,k}$ the \emph{Skorokhod (or the divergence-type) integral} of $g$ with respect to the HP $Z^{H,k}$. If we need to consider the integral as a process, we define $\int_0^\bullet g_s\delta Z_s^{H,k} = (\int_0^t g_s\delta Z_s^{H,k})_{t\in [0,T]}$ by 
	\[
		 \int_0^t g_s\delta Z_s^{H,k} := \int_0^T \bm{1}_{[0,t)}(s)g_s\delta Z_s^{H,k}, \quad t\in [0,T].
	\]

\section{$1/H$-variation of the stochastic integral with respect to HP}
\label{sec:main}

In this section, we prove the main result of the article. To this end, we denote 
	\[ V_{n,T}^p(X):= \sum_{i=0}^{2^n-1} |X_{t_{i+1}^n} -X_{t_i^n}|^p\]
for $p\in (0,\infty)$, $n\in\N$, the dyadic partition $\{t_i^n:=Ti2^{-n}\}_{i=0}^{2^n}$ of the interval $[0,T]$, and a stochastic process $X=(X_t)_{t\in [0,T]}$. Let also $H\in (1/2,1)$ and $k\in\N$ be fixed for the rest of the section. The main result of the article is stated now.

\begin{theorem}
\label{thm:main}
Let $g\in\L^{k,\frac{1}{H}}$. There is the convergence
	\[ V_{n,T}^\frac{1}{H}\left(\int_0^\bullet g_s\delta Z_s^{H,k}\right) \quad\xrightarrow[n\to\infty]{L^1(\Omega)} \quad C_{H,k} \int_0^T |g_s|^\frac{1}{H}\d{s}\]
where $C_{H,k} := \E|Z_1^{H,k}|^\frac{1}{H}$. 
\end{theorem}

In order to prove \autoref{thm:main}, we start by a slight strengthening of \cite[Proposition 6]{AssTud20}; namely we show that the convergence $V_{n,T}(Z^{H,k}) \to C_{H,k} T $ is not only in probability but also in $L^1(\Omega)$. 

\begin{lemma}
\label{lem:var_of_RP}
There is the convergence
	\[ V_{n,T}^\frac{1}{H}(Z^{H,k}) \quad\xrightarrow[n\to\infty]{L^1(\Omega)} \quad C_{H,k}T.\]
\end{lemma}

\begin{proof}
By \cite[Proposition 6]{AssTud20}, there is the following convergence in probability:
	\[ V_{n,T}^\frac{1}{H}(Z^{H,k}) \quad\xrightarrow[n\to\infty]{\P} \quad C_{H,k}T.\]
Therefore, in order to show that the convergence is also in $L^1(\Omega)$, it suffices to show that the sequence $\{V_{n,T}^\frac{1}{H}(Z^{H,k})\}_{n\in\N}$ is uniformly integrable. To this end, note that using the generalized Minkowski inequality, self-similarity, and stationarity of the increments of $Z^{H,k}$ successively yields, for any $p\in [1,\infty)$, the estimate 
\begin{align*}
	\left(\E \left(V_{n,T}^\frac{1}{H}(Z^{H,k})\right)^p \right)^\frac{1}{p}  
    & = \left(\E \left( \sum_{i=0}^{2^n-1} \left|Z^{H,k}_{t^n_{i+1}}-Z^{H,k}_{t^n_{i}}\right|^\frac{1}{H} \right)^p\right)^\frac{1}{p} \\
    & \leq \sum_{i=0}^{2^n-1}\left( \E  \left|Z^{H,k}_{t^n_{i+1}}-Z^{H,k}_{t^n_{i}}\right|^\frac{p}{H}\right)^\frac{1}{p} \\
    &= \sum_{i=0}^{2^n-1}\left( \E  \left|Z^{H,k}_{\frac{(i+1)T}{2^n}}-Z^{H,k}_{\frac{iT}{2^n}}\right|^\frac{p}{H}\right)^\frac{1}{p} \\
    &= \sum_{i=0}^{2^n-1}\left( \E  \left(T^H 2^{-nH}\left|Z^{H,k}_{i+1}-Z^{H,k}_{i}\right|\right)^\frac{p}{H}  \right)^\frac{1}{p} \\
   &= \sum_{i=0}^{2^n-1}T\, 2^{-n}\left( \E  |Z^{H,k}_1|^\frac{p}{H}  \right)^\frac{1}{p}\\
    &=T \left(\E |Z^{H,k}_1|^\frac{p}{H}\right)^\frac{1}{p}.
\end{align*}
As the last expression is finite, uniform integrability of the sequence $\{V_{n,T}^\frac{1}{H}(Z^{H,k})\}_{n\in\N}$ follows.
\end{proof}

The proof of \autoref{thm:main} can be given now.

\begin{proof}[Proof of \autoref{thm:main}]
Since, by \autoref{lem:density_of_E}, the space $\mathcal{E}_{\infty,b}$ is dense in $\L^{k,\frac{1}{H}}$, we can find a sequence of elementary processes $\{g^m\}_{m} \subseteq \mathcal{E}_{\infty,b}$ of the form
     \begin{equation*}
        g^m = \sum_{j=0}^{m-1} F_j^m\bm{1}_{[s_j^m,s_{j+1}^m)},
    \end{equation*}
    where $\{F_j^m\}_{j=0}^{m-1} \subset \mathcal{S}_{\infty,b}$ and $\{s_j^m\}_{j=0}^{m}$ is a partition of interval $[0,T]$, such that
    \begin{equation}
    \label{eq:approx}
        \|g-g^m\|_{\L^{k,\frac{1}{H}}} \xrightarrow[m \to \infty]{} 0.
    \end{equation}
Let $\{X^m\}_{m\in\N}$ be the sequence of stochastic processes $X^m=(X_t^m)_{t\in [0,T]}$, $m\in\N$, defined by $X^m_t := \int_0^t g_s^m \delta{Z^{H,k}_s}$, $t\in [0,T]$. Then we have
    \begin{align*}
        \E\left|V_{n,T}^{\frac{1}{H}}(X) - C_{H,k} \int_0^T |g_s|^\frac{1}{H}\d{s}\right| &\leq \E\left|V_{n,T}^{\frac{1}{H}}(X) - V_{n,T}^{\frac{1}{H}}(X^m)\right| \\
            &\quad + \E\left|V_{n,T}^{\frac{1}{H}}(X^m) - C_{H,k} \int_0^T |g_s^m|^\frac{1}{H}\d{s}\right| \\
            &\quad + C_{H,k}\E\left|\int_0^T (|g_s^m|^\frac{1}{H}-|g_s|^\frac{1}{H})\d{s}\right| \\
            &=: a_n^m + b_n^m + c^m
    \end{align*}
    for $m,n\in\N$ by the triangle inequality and we need to show that all the three terms in the expression above converge to zero. We will prove this in multiple steps.
    
    \textit{Step 1.}\quad We show that $\sup_{n\in\N} a_n^m \xrightarrow[m\to\infty]{} 0$ first. By \autoref{lem: est var int}, we have
    \[
     \sup_{n \in \N}   a_n^m = \sup_{n \in \N} \E\left|V_{n,T}^{\frac{1}{H}}(X) - V_{n,T}^{\frac{1}{H}}(X^m)\right|
     \lesssim \|g - g^m\|_{\L^{k,\frac{1}{H}}}\left(\|g\|_{\L^{k,\frac{1}{H}}}^\frac{1-H}{H}+ \|g^m\|_{\L^{k,\frac{1}{H}}}^\frac{1-H}{H}\right)
    \]
which converges to zero as $m\to\infty$ by \eqref{eq:approx} and by the fact that the sequence $\{\|g^m\|_{\L^{k,\frac{1}{H}}}\}_m$ is convergent and, therefore, bounded.
  
\textit{Step 2.}\quad Now we show that $c^m\xrightarrow[m\to\infty]{}0$.  By applying the Mean Value Theorem and using the H\"{o}lder inequality, we obtain the estimate
    \begin{align*}
        c^m &= C_{H,k}\E\left|\int_0^T (|g_s^m|^\frac{1}{H}-|g_s|^\frac{1}{H})\d{s}\right|\\
        & \lesssim \E \int_0^T |g_s^m-g_s|(|g_s^m|^{\frac{1}{H}-1}+|g_s|^{\frac{1}{H}-1}) \d{s}\\
        & \leq \left(\E \int_0^T |g_s^m - g_s|^\frac{1}{H}\d{s}\right)^H\left(\E\int_0^T |g_s^m|^\frac{1}{H}\d{s}\right)^{1-H} \\
        & \qquad + \left(\E \int_0^T |g_s^m - g_s|^\frac{1}{H}\d{s}\right)^H\left(\E\int_0^T |g_s|^\frac{1}{H}\d{s}\right)^{1-H}
    \\
    & \leq \|g^m-g\|_{\L^{k,\frac{1}{H}}}(\|g^m\|_{\L^{k,\frac{1}{H}}}^\frac{1-H}{H}+\|g\|_{\L^{k,\frac{1}{H}}}^\frac{1-H}{H})
    \end{align*}
where we used the embedding $\D^{k,\frac{1}{H}} \hookrightarrow L^\frac{1}{H}(\Omega)$. The right-hand side of the above inequality again converges to zero as $m\to\infty$ by the same argument as in \emph{Step 1}.
    
\textit{Step 3.}\quad Finally, we show there is the convergence $b_n^m\xrightarrow[n\to\infty]{} 0$ for every $m\in\N$. To this end, let $m\in\N$ be fixed and let us assume, without loss of generality, that $\{s_j^m\}_j\subseteq \{t_{i}^n\}_i$ for sufficiently large $n\in\N$. Denote 
	\[
		j(i):=\max \{j\in \{0,\ldots, m-1\}\mid s_j^m\leq t_i^n\}, \quad i\in\{0,\ldots,2^n\},
	\]
and 
	\[ 
		I_n^m(j):=\{i\in\{0,\ldots, 2^n\}\mid t_{i}^n\in [s_j^m, s_{j+1}^m)\}, \quad j\in \{0, \ldots, m-1\}.
	\]
Note that if $\varXi=(\varXi_t)_{t\in [0,T]}$ is a stochastic process defined by 
	\[ \varXi_t = \sum_{j=0}^{m-1} F_j^m (\xi_{s_{j+1}^m\wedge t}-\xi_{s_{j}^m\wedge t}), \quad t\in [0,T],\]
for some other stochastic process $\xi=(\xi_t)_{t\in [0,T]}$, we have that 
	\begin{equation}
	\label{eq:var_of_step}
		V_{n,T}^\frac{1}{H}(\varXi) 
			= \sum_{i=0}^{2^n-1} |F_{j(i)}^m|^\frac{1}{H} |\xi_{t_{i+1}^n}-\xi_{t_i^n}|^\frac{1}{H} 
			= \sum_{j=0}^{m-1} |F_{j}^m|^\frac{1}{H} \sum_{i\in I_n^m(j)} |\xi_{t_{i+1}^n}-\xi_{t_i^n}|^\frac{1}{H}.
	\end{equation}

\textit{Step 3.1}\quad We prove the claim for deterministic constants $\{F_j^m\}_{j}$ first. By \eqref{eq:var_of_step}, we obtain 
	\begin{align*}
		b_n^m 
			& = \E \left| \sum_{j=0}^{m-1} |F_j^m|^\frac{1}{H} \left( \sum_{i\in I_n^m(j)} |Z_{t_{i+1}^n}^{H,k} - Z_{t_{i}^n}^{H,k}|^\frac{1}{H} - C_{H,k} (s_{j+1}^m-s_j^m)\right)\right| \\
			& = \E \left| \sum_{j=0}^{m-1} |F_j^m|^\frac{1}{H} \left( \sum_{i\in I_n^m(j)} |Z_{t_{i+1}^n}^{H,k} - Z_{t_{i}^n}^{H,k}|^\frac{1}{H} - \E|Z_{s_{j+1}^m}^{H,k} - Z_{s_{j}^m}^{H,k}|^\frac{1}{H}\right)\right|\\
			& \leq \sum_{j=0}^{m-1} |F_j^m|^\frac{1}{H} \E\left| \sum_{i\in I_n^m(j)} |Z_{t_{i+1}^n}^{H,k} - Z_{t_{i}^n}^{H,k}|^\frac{1}{H} - \E|Z_{s_{j+1}^m}^{H,k} - Z_{s_{j}^m}^{H,k}|^\frac{1}{H}\right|.
	\end{align*}
Since we have the convergence
	\[ 
		\E\left| \sum_{i\in I_{n}^m(j)} |Z_{t_{i+1}^n}^{H,k} - Z_{t_{i}^n}^{H,k}|^\frac{1}{H} - \E|Z_{s_{j+1}^m}^{H,k} - Z_{s_{j}^m}^{H,k}|^\frac{1}{H}\right|\xrightarrow[n\to\infty]{} 0
	\]
for every $j\in\{0,\ldots, m-1\}$ by \autoref{lem:var_of_RP}, self-similarity of the HP $Z^{H,k}$, and stationarity of its increments, the convergence $b_n^m\xrightarrow[n\to\infty]{} 0$ follows. This, together with \textit{Step 1} and \textit{Step 2}, concludes the proof of \autoref{thm:main} for deterministic functions $g\in L^{\frac{1}{H}}([0,T])$.

\textit{Step 3.2}\quad Let us now consider the general case of $\{F_j^m\}_j\subset \mathcal{S}_{\infty,b}$. By \autoref{lem:expansion}, we have 
	\begin{equation}
	\label{eq:expansion_used}
		X_t^m = \sum_{\ell=0}^k (-1)^\ell{k\choose \ell} \sum_{j=0}^{m-1}\int_{[0,T]^\ell}(D^\ell F_j^m)(x)\delta^{k-\ell} ((L^{H,k}\bm{1}_{[s_j^m\wedge t, s_{j+1}^m\wedge t)})(x,\bullet)) \d{x}, \quad t\in [0,T].
	\end{equation}
We first note that for $\ell\in\{0,\ldots, k-1\}$, $j\in \{0,\ldots, m-1\}$, and almost every $x\in [0,T]^\ell$, we have 
	\[
		\delta^{k-\ell} ((L^{H,k}\bm{1}_{[s_j^m\wedge t, s_{j+1}^m\wedge t)})(x,\bullet)) \\ =  \int_{s_j^m\wedge t}^{s_{j+1}^m\wedge t} g^\ell(x,u)\delta Z_u^{H',k-\ell}, \quad t\in [0,T].
	\]
where 
	\[
		H':=H\left(1-\frac{\ell}{k}\right) + \frac{\ell}{k}
	\]
and
	\[ g^\ell (y_1,\ldots,y_\ell,v) := c_\ell \prod_{i=1}^\ell \left(\frac{v}{y_i}\right)^{\frac{1}{2}-\frac{1-H}{k}} (v-y_i)_+^{-\left(\frac{1}{2}+\frac{1-H}{k}\right)}
	\]
for $y_1,\ldots, y_\ell\in (0,T]$ and $v\in [0,T]$, in which 
	\[
		 c_\ell := \frac{c_{H,k}}{c_{H',k-\ell}}.
	\] 
This means, that if we write $X^m =: \sum_{\ell=0}^k Y^{m,\ell}$, each term $Y^{m,\ell}$, $\ell \in\{0,\ldots,k\}$, in the expansion \eqref{eq:expansion_used} can be thought of as a mixed Lebesgue--Skorokhod integral with respect to the HP of order $k-\ell$ with Hurst index $H'$. In particular, we have 
	\begin{align*}
		Y_t^{m,0} & =  \sum_{j=0}^{m-1} F_j^m(Z_{s_{j+1}^m\wedge t}^{H,k} - Z_{s_{j}^m\wedge t}^{H,k}), \\
		Y_t^{m,1} & =  -kc_1 \sum_{j=0}^{m-1}\int_0^T (DF^m_j)(x) \left(\int_{s_j^m\wedge t}^{s_{j+1}^m\wedge t} \left(\frac{u}{x}\right)^{\frac{1}{2}-\frac{1-H}{k}} (u-x)_+^{-\left(\frac{1}{2} + \frac{1-H}{k}\right)}\delta Z_u^{H\left(1-\frac{1}{k}\right)+\frac{1}{k},k-1}\right)\d{x},\\
		& \vdots  \\
		Y_t^{m,k} & =  (-1)^k c_k\sum_{j=0}^{m-1} \int_{[0,T]^k} (D^kF_j^m)(x) \left(\int_{s_j^m\wedge t}^{s_{j+1}^m\wedge t} \prod_{i = 1}^k \left(\frac{u}{x_i}\right)^{\frac{1}{2}-\frac{1-H}{k}} (u-x_i)_+^{-\left(\frac{1}{2}+\frac{1-H}{k}\right)}\d{u}\right)\d{x},
	\end{align*}
where we set, for convenience, $c_k:=c_{H,k}$. We aim to show that the $1/H$-variation of $X^m$ is governed by the $1/H$-variation of the term $Y^{m,0}$. 
Let us write 
	\[ b_n^m \leq \E\left| V_{n,T}^\frac{1}{H}(X^m) - V_{n,T}^\frac{1}{H}(Y^{m,0})\right| + \E\left|V_{n,T}^\frac{1}{H}(Y^{m,0}) - C_{H,k}\int_0^T|g^m|^\frac{1}{H} \d{s}\right|=: d_n^m + e_n^m.\]

	\textit{Step 3.2.1}\quad We prove that $d^m_n\xrightarrow[n\to\infty]{} 0$ first. Define $R^m=(R_t^m)_{t\in [0,T]}$ by
	\[
		R^m_t := \sum_{\ell=1}^{k} Y_t^{m,\ell},\quad t\in [0,T].
	\]
By using \autoref{lem: est var}, we obtain 
		\[ 
			d_n^m \lesssim (\E V_{n,T}^\frac{1}{H}(R^m))^H ( (\E V_{n,T}^\frac{1}{H}(X^m))^{1-H} + (\E V_{n,T}^\frac{1}{H}(Y^{m,0}))^{1-H}).
		\]
	By using \autoref{lem: est var int}, we have the estimate 
		\[ 
		  \E V_{n,T}^\frac{1}{H}(X^m) \lesssim \|g^m\|_{\L^{k,\frac{1}{H}}}^\frac{1}{H}.
		\]
By \eqref{eq:var_of_step} and \autoref{lem: est var int}, we also have 
		\begin{align*}
			\E V_{n,T}^\frac{1}{H}(Y^{m,0}) 
				& = \E \sum_{i=0}^{2^n-1} \left|F_{j(i)}^m\right|^\frac{1}{H}\left|Z_{t_{i+1}^n}^{H,k}-Z_{t_i^n}^{H,k}\right|^\frac{1}{H} \\
				& \leq \max_{j\in\{0,\ldots,m-1\}} \sup_{\omega\in\Omega}|F_j^m(\omega)|^\frac{1}{H} \E V_{n,T}^\frac{1}{H}(Z^{H,k}) \\
				& \lesssim \max_{j\in\{0,\ldots,m-1\}} \sup_{\omega\in\Omega}|F_j^m(\omega)|^\frac{1}{H} T.
		\end{align*}
It remains to show that 
	\begin{equation}
	\label{eq:VnRm}
		\E V_{n,T}^\frac{1}{H}(R^m) \xrightarrow[n\to\infty]{} 0.
	\end{equation}
By \autoref{lem:triangle_for_var}, we have that 
	\[ 
		V_{n,T}^\frac{1}{H}(R^m) \lesssim \sum_{\ell=1}^k V_{n,T}^\frac{1}{H}(Y^{m,\ell})
	\]
so that the convergence in \eqref{eq:VnRm} will follow, once we show that 
	\begin{equation}
	\label{eq:varYml}
	\E V_{n,T}^\frac{1}{H}(Y^{m,\ell})\xrightarrow[n\to\infty]{} 0
	\end{equation}
holds for every $\ell\in \{1,\ldots,k\}$. We first prove convergence \eqref{eq:varYml} for $\ell\in \{1,\ldots, k-1\}$ and then for $\ell=k$ separately. 

Let $\ell\in \{1,\ldots, k-1\}$. We have, by using the Minkowski inequality, 
	\begin{align*}
		(\E V_{n,T}^\frac{1}{H}(Y^{m,\ell}))^H & \propto \left(\E \sum_{i=0}^{2^n-1} \left| \int_{[0,T]^\ell} (D^\ell F_{j(i)}^m)(x) \int_{t_{i}^n}^{t_{i+1}^n} g^\ell(x,u) \delta Z_u^{H',k-\ell} \d{x}\right|^\frac{1}{H} \right)^H \\
			& \lesssim \int_{[0,T]^\ell} \left(\E\sum_{i=0}^{2^n-1} \left| (D^\ell F_{j(i)}^m)(x)\int_{t_{i}^n}^{t_{i+1}^n} g^\ell (x,u)\delta Z_u^{H',k-\ell} \right|^\frac{1}{H} \right)^H\d{x} \\
			& \leq \sup_{j\in\{0,\ldots, m-1\}} \sup_{\omega\in\Omega}\sup_{x\in [0,T]^\ell} |(D^\ell F_j^m)(\omega,x)| \\
			& \qquad  \times\int_{[0,T]^\ell} \left(\E V_{n,T}^\frac{1}{H}\left(\int_0^\bullet g^\ell(x,u)\delta Z_u^{H',k-\ell} \right)\right)^H \d{x}.\numberthis\label{eq:bound}
	\end{align*}
Now, it is easily verified that, for almost every $x\in [0,T]^\ell$, the function $g^\ell(x,\bullet)$ belongs to $L^\frac{1}{H'}([0,T])$ so that the process
	\[ 
		U(x) := \int_0^\bullet g^\ell(x,u)\delta Z_u^{H', k-\ell}
	\]
has finite $1/H'$-variation as we have already proved \autoref{thm:main} for deterministic integrands (see the end of \emph{Step 3.1}), i.e. we have the convergence
	\begin{equation}
	\label{eq:1/H'_var_of_U}
		\E V_{n,T}^\frac{1}{H'}(U(x)) \xrightarrow[n\to\infty]{} C_{H',k-\ell} \int_0^T \left|g^{\ell}(x,u)\right|^\frac{1}{H'}\d{u}
	\end{equation}
Moreover, we have that $H'>H$ and, by the equivalence of moments of random variables on a finite Wiener chaos (cf., e.g., \cite[Corollary 2.8.14]{NouPec12}), we obtain the estimate
    \begin{align*}
        \E V_{n,T}^\frac{1}{H}(U(x)) 
        	& =  \sum_{i=0}^{2^n-1} \E \left|U_{t^n_{i+1}}(x)-U_{t^n_{i}}(x)\right|^\frac{1}{H} \\
        	& \eqsim \sum_{i=0}^{2^n-1} \left(\E \left|U_{t^n_{i+1}}(x)-U_{t^n_{i}}(x)\right|^\frac{1}{H'}\right)^\frac{H'}{H} \\
       		& \leq \sup_{i=0,\ldots, 2^n-1}\left(\E \left|U_{t^n_{i+1}}(x)-U_{t^n_{i}}(x)\right|^\frac{1}{H'}\right)^{\frac{H'}{H}-1} \left( \sum_{i=0}^{2^n-1}\E \left|U_{t^n_{i+1}}(x)-U_{t^n_{i}}(x)\right|^\frac{1}{H'} \right)\\
        	& = \left(\sup_{i=0,\ldots, 2^n-1} \E \left|U_{t^n_{i+1}}(x)-U_{t^n_{i}}(x)\right|^\frac{1}{H'}\right)^{\frac{H'}{H}-1} \left( \E V_{n,T}^\frac{1}{H'}(U(x)) \right) \numberthis\label{eq:1/H_var_of_U}
    \end{align*}
To analyze the first factor in \eqref{eq:1/H_var_of_U}, we can apply \autoref{lem:bddness_of_integral} to obtain
    \[
        \E \left|U_{t^n_{i+1}}(x)-U_{t^n_{i}}(x)\right|^\frac{1}{H'} =\E \left| \int_{t^n_{i}}^{t^n_{i+1}} g^\ell(x,u)\delta Z_u^{H', k-\ell}\right|^\frac{1}{H'}  \lesssim \int_{t^n_{i}}^{t^n_{i+1}} \left| g^\ell(x,u)\right|^\frac{1}{H'}\d u 
    \]
and integrability of $\left|g^\ell(x,\bullet)\right|^\frac{1}{H'}$ and absolute continuity of (Lebesque) integral guarantee that
    \[
        \sup_{i=0,\ldots 2^n-1} \left( \int_{t^n_{i}}^{t^n_{i+1}} \left| g^\ell(x,u)\right|^\frac{1}{H'}\d u \right)\xrightarrow[n\to\infty]{} 0.
    \]
The second factor in \eqref{eq:1/H_var_of_U} converges by \eqref{eq:1/H'_var_of_U}. Altogether, it follows that 
	\[ 
		\E V_{n,T}^{\frac{1}{H}}(U(x)) \xrightarrow[n\to\infty]{} 0.
	\]
Now, by \autoref{lem:bddness_of_integral}, we obtain the uniform bound
	\begin{align*}
		\E V_{n,T}^\frac{1}{H} (U(x))  
			&  =  \sum_{i=0}^{2^n-1} \E\left| \int_{t_i^n}^{t_{i+1}^n} g^\ell(x,u)\delta Z_u^{H', k-\ell}\right|^\frac{1}{H} \\
			&  \lesssim \sum_{i=0}^{2^n-1}\left(\int_{t_i^n}^{t_{i+1}^n} |g^\ell (x,u)|^\frac{1}{H'} \d{u}\right)^{\frac{H'}{H}}\\
			&  = \sum_{i=0}^{2^n-1} \left(\int_{t_i^n}^{t_{i+1}^n} |g^\ell (x,u)|^\frac{1}{H'} \d{u}\right)^{\frac{H'}{H}-1}\left(\int_{t_i^n}^{t_{i+1}^n} |g^\ell (x,u)|^\frac{1}{H'} \d{u}\right)\\
			&  \leq \left(\int_{0}^{T} |g^\ell (x,u)|^\frac{1}{H'} \d{u}\right)^{\frac{H'}{H}-1} \sum_{i=0}^{2^n-1}\int_{t_i^n}^{t_{i+1}^n} |g^\ell (x,u)|^\frac{1}{H'} \d{u} \\
			&  = \left(\int_{0}^{T} |g^\ell (x,u)|^\frac{1}{H'} \d{u}\right)^{\frac{H'}{H}}
	\end{align*} 
and as we have that
	\[
		\int_{[0,T]^\ell} \left(\int_{0}^{T} |g^\ell (x,u)|^\frac{1}{H'} \d{u}\right)^{H'}\d{x}<\infty
	\]
holds by \autoref{lem:convergence of integral g}, convergence \eqref{eq:varYml} for $\ell\in \{1,2,\ldots, k-1\}$ follows by taking the limit $n\to\infty$ of  \eqref{eq:bound} and appealing to the Dominated Convergence Theorem. 

Convergence \eqref{eq:varYml} for $\ell=k$ is proved now. We begin with the estimate
	\[
		\E V_{n,T}^\frac{1}{H}(Y^{m,k}) \leq \sup_{j\in\{0,\ldots, m-1\}}\sup_{\omega\in\Omega}\sup_{x\in[0,T]^\ell} |(D^kF_j^m)(\omega,x)|^\frac{1}{H}\sum_{i=0}^{2^n-1} \left(\int_{[0,T]^k}\int_{t_i^n}^{t_{i+1}^n} |g^k(x,u)|\d{x}\d{u}\right)^\frac{1}{H}.
	\]
Now, we see that 
	\[ \left(\int_{[0,T]^k} \int_0^T g^k(x,u) \d{u}\d{x}\right)^2 \lesssim \int_{[0,T]^k}\left( \int_0^T g^k(x,u) \d{u}\right)^2\d{x} = \int_{[0,T]^k} |( L^{H,k} 1)(x) |^2\d{x}<\infty\]
holds by the Cauchy--Schwarz inequality and \autoref{lem:bddness_of_LH}. It follows that the function
	\[ 
		\int_{[0,T]^k}\int_0^\bullet g^k(x,u)\d{u}\d{x} = \int_0^\bullet \int_{[0,T]^k} g^k(x,u)\d{x}\d{u},
	\]
being increasing and integrable, is of bounded variation and, therefore, of zero $1/H$-variation. Thus we obtain 
	\[
		\sum_{i=0}^{2^n-1} \left(\int_{[0,T]^k}\int_{t_i^n}^{t_{i+1}^n} |g^k(x,u)|\d{x}\d{u}\right)^\frac{1}{H} \xrightarrow[n\to\infty]{} 0
	\] 
and convergence \eqref{eq:varYml} for $\ell=k$, and therefore also convergence $d_n^m\xrightarrow[n\to\infty]{}0$, is proved. 

\textit{Step 3.2.2}\quad As in \emph{Step 3.1}, we obtain the estimate 
	\[
		e_n^m \leq \sum_{j=0}^{m-1} \sup_{\omega\in\Omega} |F_j^m(\omega)|^\frac{1}{H} \E\left| \sum_{i\in I_n^m(j)} |Z_{t_{i+1}^n}^{H,k} - Z_{t_{i}^n}^{H,k}|^\frac{1}{H} - \E|Z_{s_{j+1}^m}^{H,k} - Z_{s_{j}^m}^{H,k}|^\frac{1}{H}\right|
	\]
and the convergence $e_n^m\xrightarrow[n\to\infty]{} 0$ now follows by the same arguments used there.
\end{proof}

\section*{Declarations}

\noindent \textbf{Declaration of interest:} None.

\noindent \textbf{Funding:} Petr Čoupek was supported by the Czech Science Foundation project no.\ 23-05737S. Pavel Kříž was supported by the Czech Science Foundation project no.\ 24-11146S.

\appendix

\section{Auxiliary lemmas}
\label{sec:app}

\subsection{An auxiliary lemma on boundedness of the divergence-type integral}

Recall that we denote by $\mathcal{S}_{\infty,b}$ the space of random variables $F$ such that $F=f(I(h_1),\ldots, I(h_n))$ for some $n\in\N$, $f\in C_b^\infty(\R^n)$, and $h_1,\ldots, h_n\in C([0,T])$. Recall also that we denote by $\mathcal{E}_{\infty,b}$ the space of $\mathcal{S}_{\infty,b}$-valued elementary processes defined on the interval $[0,T]$, i.e.\ the space of all processes $g$ of the form 
	\begin{equation}
	\label{eq:elem_s}
		g = \sum_{i=0}^{n-1} G_i\bm{1}_{[s_i, s_{i+1})},
	\end{equation}
where $n\in\N$, $\{G_i\}_{i=0}^{n-1} \subset \mathcal{S}_{\infty,b}$, and where $\{0=s_0\leq s_1\leq \ldots \leq s_n=T\}$ is a partition of the interval $[0,T]$. The following lemma is used for the construction of the divergence-type integral with respect to an HP.

\begin{lemma}
\label{lem:boundedness_of_int_appendix}
Let $H\in (1/2,1)$ and $K\in\N$ and let $Z^{H,K}$ denote the HP of order $K$ with Hurst parameter $H$. Then for every $k\in\N$ such that $k\geq K$, every $p,q$ such that $1\leq pH<\infty$ and $1\leq qH<\infty$, for every $0\leq u\leq v\leq T$, and for every $g\in\mathcal{E}_{\infty,b}$, there is the estimate
	\begin{equation*}
	\left\|\int_u^v g_s\delta Z_s^{H,K}\right\|_{\D^{k-K,q}}\lesssim (v-u)^{H-\frac{1}{p}}\|g\|_{L^{p}([u,v],\D^{k,q})}.
	\end{equation*}
\end{lemma}

\begin{proof}
The proof is similar to that of \cite[Proposition 3]{CouDunDun22}. By using the Meyer inequalities from, e.g., \cite[Theorem 2.5.5]{NouPec12} or \cite[Proposition 1.5.7]{Nua06}, we obtain the estimate 
	\[
		\left\| \int_0^T g_s\delta Z_s^{H,K}\right\|_{\D^{k-K,q}}^q \lesssim \sum_{\ell=0}^k \E \|(D^\ell L^{H,K}g({\color{blue} \bullet}))({\color{red}\bullet}) \|_{L^2({\color{blue}[0,T]^K})\otimes L^2({\color{red}[0,T]^\ell})}^q.
	\]
For every $\ell\in \{0,\ldots, k\}$, we have
	\begin{align*}
		\E \|(D^\ell L^{H,K}g({\color{blue} \bullet}))({\color{red}\bullet}) \|_{L^2({\color{blue}[0,T]^K})\otimes L^2({\color{red}[0,T]^\ell})}^q & \\
			& \hspace{-4cm} = \E\left(\int_{[0,T]^\ell} \left\| L^{H,K}D^\ell g(x) \right\|_{L^2([0,T]^K)}^2 \d{x}\right)^{\frac{q}{2}}\\ 
			& \hspace{-4cm} \lesssim \E \left(\int_{[0,T]^\ell} \|D^\ell g(x)\|_{L^{\frac{1}{H}}([0,T])}^{2}\d{x}\right)^\frac{q}{2} \\
			& \hspace{-4cm} = \E \left(\int_{[0,T]^\ell} \left(\int_0^T |(D^\ell g_s)(x)|^\frac{1}{H}\d{s}\right)^{2H}\d{x}\right)^\frac{q}{2}\\
			& \hspace{-4cm} \leq \E \left(\int_0^T \left(\int_{[0,T]^\ell} |(D^\ell g_s)(x)|^2\d{x}\right)^\frac{1}{2H} \d{s}\right)^{qH}\\
			& \hspace{-4cm} \leq \left( \int_0^T \left(\E \left( \int_{[0,T]^\ell} |(D^\ell g_s)(x)|^2\d{x}\right)^\frac{q}{2}\right)^\frac{1}{qH}\d{s}\right)^{qH} \\
			& \hspace{-4cm} \leq \left(\int_0^T \|g_s\|_{\D^{k,q}}^\frac{1}{H}\d{s}\right)^{qH}
	\end{align*}
by using \autoref{lem:bddness_of_LH} and the Minkowski inequality twice successively. Thus we obtain the estimate
	\[ 
		\left\|\int_0^T g_s\delta Z_s^{H,K}\right\|_{\D^{k-K,q}} \lesssim \|g\|_{L^\frac{1}{H}([0,T];\D^{k,q})}.
	\]
The claim is then proved by applying this inequality to the process $g\bm{1}_{[u,v)}\in \mathcal{E}_{\infty,b}$ and subsequently using the H\"older inequality.
\end{proof}

\subsection{Auxiliary lemmas on Malliavin calculus}

The following lemma is a density result that is used for the construction of the divergence-type integral with respect to an HP but also in the proof of the main \autoref{thm:main}. 

\begin{lemma}
\label{lem:density_of_E_appendix}
For every $k\in\N$ and $p,q\in [1,\infty)$, the space $\mathcal{E}_{\infty,b}$ is dense in the space $L^q([0,T],\D^{k,p})$.
\end{lemma}

\begin{proof}

We split the proof into two steps. 

\emph{Step 1.}\quad First we show that the space $\mathcal{S}_{\infty,b}$ is dense in $\D^{k,p}$. To this end, recall that the space $\mathcal{S}_b$, i.e.\ the space or random variables $F$ such that $F=f(I(h_1),\ldots, I(h_n))$ for some $n\in\N$, $f\in C_b^\infty(\R^n)$, and $h_1,\ldots, h_n\in L^2([0,T])$, is dense in $\D^{k,p}$ (cf.\ e.g., \cite[p. 25]{Nua06}) and so the claim would follow if we show that each element $F\in \mathcal{S}_b$ can be approximated by a sequence of elements from space $\mathcal{S}_{\infty,b}$. 

Let $F\in\mathcal{S}_b$ be of the form $F=f(I(h_1),\ldots, I(h_n))$ for some $n\in\N$, $f\in C_b^\infty(\R^n)$, and $h_1,\ldots, h_n\in L^2([0,T])$. Since continuous functions are dense in $L^2([0,T])$, for each $i\in \{1,\ldots, n\}$, we can find a sequence $\{h_i^m\}_m$ of continuous functions such that $\| h_i^m-h_i\|_{L^2([0,T])} \xrightarrow[m\to\infty]{} 0$. Define $F^m := f(I(h_1^m),\ldots, I(h_n^m))$. We clearly have that $F^m\in\mathcal{S}_{\infty,b}$ and we aim to show the convergence 
	\begin{equation}
	\label{eq:density_0}
		\|F-F^m\|_{\D^{k,p}}\xrightarrow[m\to\infty]{} 0.
	\end{equation}
To this end, note that we have, for every $\ell\in\{0,\ldots, k\}$, the estimate 
	\begin{align*}
		\E\|D^\ell F-D^\ell F^m\|_{L^2([0,T]^\ell)}^p & \\
			& \hspace{-3cm} = \E\Bigg\| \sum_{i_1=1}^n\cdots\sum_{i_\ell=1}^n\Big[ (\partial_{i_1,\ldots, i_\ell}^\ell f)(I(h_1),\ldots, I(h_n))h_{i_1}\otimes \ldots \otimes h_{i_\ell} \\
			& \hspace{0cm} - (\partial_{i_1,\ldots, i_\ell}^\ell f)(I(h_1^m),\ldots, I(h_n^m))h_{i_1}^m\otimes \ldots \otimes h_{i_\ell}^m \Big]\Bigg\|_{L^2([0,T]^\ell)}^p \\
			& \hspace{-3cm} \lesssim \sum_{i_1=1}^n\cdots\sum_{i_\ell=1}^n \Big[ \E\Big\| (\partial_{i_1,\ldots, i_\ell}^\ell f)(I(h_1),\ldots, I(h_n))h_{i_1}\otimes \ldots \otimes h_{i_\ell} \\
			& \hspace{0cm} - (\partial_{i_1,\ldots, i_\ell}^\ell f)(I(h_1),\ldots, I(h_n))h_{i_1}^m\otimes \ldots \otimes h_{i_\ell}^m\Big\|_{L^2([0,T]^\ell)}^p \\
			& + \E \Big\| (\partial_{i_1,\ldots, i_\ell}^\ell f)(I(h_1),\ldots, I(h_n))h_{i_1}^m\otimes \ldots \otimes h_{i_\ell}^m  \\
			& \hspace{0cm} - (\partial_{i_1,\ldots, i_\ell}^\ell f)(I(h_1^m),\ldots, I(h_n^m))h_{i_1}^m\otimes \ldots \otimes h_{i_\ell}^m\Big\|_{L^2([0,T]^\ell)}^p \Big]. \numberthis\label{eq:density_1}
	\end{align*}
Now, for $(i_1, \ldots, i_\ell)\in \{1,\ldots, n\}^\ell$, we estimate the first term in \eqref{eq:density_1} by
	\begin{align*}
		 \E\Big\| (\partial_{i_1,\ldots, i_\ell}^\ell f)(I(h_1),\ldots, I(h_n))(h_{i_1}\otimes \ldots \otimes h_{i_\ell} - h_{i_1}^m\otimes \ldots \otimes h_{i_\ell}^m)\Big\|_{L^2([0,T]^\ell)}^p \\
		 & \hspace{-7cm} \leq \sup_{x\in\R^n} |(\partial_{i_1,\ldots, i_\ell}^\ell f)(x)|^p \|h_{i_1}\otimes \ldots \otimes h_{i_\ell} - h_{i_1}^m\otimes \ldots \otimes h_{i_\ell}^m\Big\|_{L^2([0,T]^\ell)}^p
	\end{align*}
and we must show that the right-hand side converges to zero as $m\to\infty$. We proceed by induction. If $\ell = 1$, the claim reduces to the claim that for every $i_1\in \{1,\ldots, n\}$, we have $\|h_{i_1}-h_{i_1}^m\|_{L^2([0,T])}\xrightarrow[m\to\infty]{} 0$ but this is exactly how $h^m_{i_1}$ is chosen. Now assume that for every $(i_1,\ldots, i_{\ell-1})\in \{1,\ldots, n\}^{\ell-1}$, there is the convergence \[\|h_{i_1}\otimes\ldots\otimes h_{i_{\ell-1}} - h_{i_1}^m\otimes\ldots\otimes h_{i_{\ell-1}}^m\|_{L^2([0,T]^{\ell-1})} \xrightarrow[m\to\infty]{} 0.\] We let $i_\ell\in\{1,\ldots, n\}$, write
	\begin{align*}
		\|h_{i_1}\otimes\ldots\otimes h_{i_{\ell}} - h_{i_1}^m\otimes\ldots\otimes h_{i_{\ell}}^m\|_{L^2([0,T]^{\ell})}
			& \\
			& \hspace{-3cm} \leq \| h_{i_1}\otimes\ldots\otimes h_{i_{\ell-1}}\otimes (h_{i_{\ell}} - h_{i_{\ell}}^m) \|_{L^2([0,T]^{\ell})} \\
			& \hspace{-2cm} + \| (h_{i_1}\otimes\ldots\otimes h_{i_{\ell-1}} - h_{i_1}^m\otimes\ldots\otimes h_{i_{\ell-1}}^m)\otimes h_{i_{\ell}}^m\|_{L^2([0,T]^{\ell})} \\
			& \hspace{-3cm} = \|h_{i_1}\otimes\ldots\otimes h_{i_{\ell-1}}\|_{L^2([0,T]^{\ell-1})} \|h_{i_{\ell}} - h_{i_{\ell}}^m\|_{L^2([0,T])} \\
			& \hspace{-2cm} + \| h_{i_1}\otimes\ldots\otimes h_{i_{\ell-1}} - h_{i_1}^m\otimes\ldots\otimes h_{i_{\ell-1}}^m \|_{L^2([0,T]^{\ell-1})} \|h_{i_\ell}^m\|_{L^2([0,T])},
	\end{align*}
and we see that both terms on the right-hand side of the above inequality tend to zero as $m\to\infty$: The convergence of the first term follows as $h_{i_\ell}^m$ were chosen so that $\|h_{i_\ell}-h_{i_{\ell}}^m\|_{L^2([0,T])}\xrightarrow[m\to\infty]{} 0$ holds and the convergence of the second term follows by using the induction hypothesis and the fact that the sequence $\{\|h_{i_\ell}^m\|_{L^2([0,T])}\}_m$ is convergent and, therefore, bounded. 

Now, for the second term in \eqref{eq:density_1}, we have that 
	\begin{align*}
		 \E | [(\partial_{i_1,\ldots, i_\ell}^\ell f)(I(h_1),\ldots, I(h_n))  - (\partial_{i_1,\ldots, i_\ell}^\ell f)(I(h_1^m),\ldots, I(h_n^m))]h_{i_1}^m\otimes \ldots \otimes h_{i_\ell}^m\Big\|_{L^2([0,T]^\ell)}^p & \\
		& \hspace{-13cm} = \E |(\partial_{i_1,\ldots, i_\ell}^\ell f)(I(h_1),\ldots, I(h_n))  - (\partial_{i_1,\ldots, i_\ell}^\ell f)(I(h_1^m),\ldots, I(h_n^m))|^p \prod_{j=1}^{\ell} \|h^m_{i_j}\|_{L^2([0,T])}^p\\
		& \hspace{-13cm} \leq \sup_{x\in\R^n} \|(\nabla \partial^\ell_{i_1,\ldots, i_{\ell}}f)(x)\|_{\R^n}^p \E \|( I(h_1-h_1^m), \ldots, I(h_n-h_n^m))\|_{\R^n}^p \prod_{j=1}^{\ell} \|h^m_{i_j}\|_{L^2([0,T])}^p.\numberthis\label{eq:density_3}
	\end{align*}
where we used the Mean Value Inequality and linearity of the Wiener--It\^o integral $I$. Furthermore, it is readily seen that
	\begin{multline*}
		 \E \|( I(h_1-h_1^m), \ldots, I(h_n-h_n^m)) \|_{\R^n}^p 
		 	 = \E \left(\sum_{i=1}^n |I(h_i-h_i^m)|\right)^{p} \\
		 	 \leq\left(\sum_{i=1}^n \left(\E|I(h_i-h_i^m)|^p\right)^\frac{1}{p}\right)^p 
		 	 \lesssim \left( \sum_{i=1}^n \left(\E|I(h_i-h_i^m)|^2\right)^\frac{1}{2}\right)^p
		 	 \lesssim \left(\sum_{i=1}^n \|h_i-h_i^m\|_{L^2([0,T])}\right)^p
	\end{multline*}
holds by using the Minkowski inequality, equivalence of moments on a finite Wiener chaos (cf., e.g., \cite[Theorem 2.7.2]{NouPec12}), and the fact that the Wiener--It\^o integral $I$ is an isometry between $L^2(\Omega)$ and $L^2([0,T])$. We thus obtain from \eqref{eq:density_3} that 
	\begin{multline}
	\label{eq:density_4}
		\E | [(\partial_{i_1,\ldots, i_\ell}^\ell f)(I(h_1),\ldots, I(h_n))  - (\partial_{i_1,\ldots, i_\ell}^\ell f)(I(h_1^m),\ldots, I(h_n^m))]h_{i_1}^m\otimes \ldots \otimes h_{i_\ell}^m\Big\|_{L^2([0,T]^\ell)}^p \\ 
		\lesssim \sup_{x\in\R^n} \|(\nabla \partial^\ell_{i_1,\ldots, i_{\ell}}f)(x)\|_{\R^n}^p \left(\sum_{i=1}^n \|h_i-h_i^m\|_{L^2([0,T])}\right)^p \prod_{j=1}^{\ell} \|h^m_{i_j}\|_{L^2([0,T])}^p.
	\end{multline}
Now, we have the convergence $\|h_i-h_i^m\|_{L^2([0,T])}\xrightarrow[m\to\infty]{}0$ for every $i\in\{1,\ldots,n\}$, the sequences $\{h_{i_j}^m\}_m$, $j\in\{1,\ldots,\ell\}$, are convergent and, therefore, bounded, and \[\sup_{x\in\R^n} \| (\nabla\partial^\ell_{i_1,\ldots, i_\ell} f)(x)\|_{\R^n}^p<\infty\] as $f\in C_b^\infty(\R^n)$. Consequently, the right-hand side of \eqref{eq:density_4} tends to zero as $m\to\infty$ which finishes the proof of \eqref{eq:density_0}. 

\emph{Step 2.}\quad  Now we show that the space $\mathcal{E}_{\infty,b}$ is dense in the space $L^q([0,T],\D^{k,p})$. To this end, denote by $\mathcal{E}^{k,p}$ the space of elementary $\D^{k,p}$-valued processes, i.e.\ the space of all processes $g$ of the form \eqref{eq:elem_s} where $\{G_i\}_{i=0}^{m-1} \subset \D^{k,p}$. It follows by standard arguments that the space $\mathcal{E}^{k,p}$ is dense in $L^q([0,T],\D^{k,p})$ and so the claim would follow if we show that each element $g\in\mathcal{E}^{k,p}$ can be approximated by a sequence of elements from $\mathcal{E}_{\infty,b}$. Let $g\in \mathcal{E}^{k,p}$ be of the form \eqref{eq:elem_s}. By step~1 of the proof, for each $G_i$, $i\in \{0,\ldots, n-1\}$, we can find a sequence $\{G_i^m\}_{m}\subset\mathcal{S}_{\infty,b}$ such that $\|G_i-G_i^m\|_{\D^{k,p}}\xrightarrow[m\to\infty]{}0$. Let us define
	\[ g^m:= \sum_{i=0}^{n-1} G_i^m\bm{1}_{[t_i,t_{i+1})}.\]
Then by straightforward computations, we obtain 
	\[ \| g-g^m\|_{L^q([0,T],\D^{k,p})}^q = \sum_{i=0}^{n-1} \|G_i-G_i^m\|_{\D^{k,p}}^q(t_{i+1}-t_i) \]
which converges to zero as $m\to\infty$.
\end{proof}

The equality in the following lemma is a generalization of \cite[Proposition 1.3.3]{Nua06} for higher-order divergence operators and it is used in the proof of the main \autoref{thm:main}. Recall that if $F\in\mathcal{S}_{\infty,b}$, then 
	\[
		\supnorm{F}{k}:=\sup_{\omega\in\Omega} \sum_{\ell=0}^k \sup_{x\in [0,T]^\ell} |(D^\ell F)(\omega,x)| <\infty
	\]
for every $k\in\N$.

\begin{lemma}
\label{lem:expansion}
Let $k\in\N$, $F\in\mathcal{S}_{\infty,b}$, and let $u\in L^2([0,T]^k)$ be a symmetric function in all variables. Then  
	\[ \delta^k(Fu) = \sum_{\ell=0}^k (-1)^\ell {k\choose \ell} \int_{[0,T]^\ell} (D^\ell F)(x) \delta^{k-\ell}(u(\bullet,x))\d{x}.\]
\end{lemma}

\begin{proof}
We proceed by induction on $k$. Note that the claim for $k=1$ follows immediately from \cite[Proposition 1.3.3]{Nua06}. Let $u\in L^2([0,T]^{k+1})$ be a symmetric function in all its ($k+1$) variables. This means that $Fu\in\Dom\delta^{k+1}$ and we also have that, for almost every $x\in [0,T]$, $u(\bullet,x)$ is a symmetric function in all its ($k$) variables that belongs to $L^2([0,T]^k)$. Therefore, we can and do assume, that 
	\[ \delta^k(Fu(\bullet,x)) = \sum_{\ell = 0}^k (-1)^\ell {k\choose \ell} \int_{[0,T]^\ell} (D^\ell F)(y)\delta^{k-\ell}(u(\bullet,y,x))\d{y}\]
holds for almost every $x\in [0,T]$. It then follows that we have 
	\begin{equation}
	\label{eq:expansion_1}
		\delta^{k+1}(Fu) = {\color{red}\delta}({\color{blue}\delta^k}(Fu({\color{blue}\bullet},{\color{red}\bullet}))) =  \sum_{\ell = 0}^k (-1)^\ell {k\choose \ell}  {\color{red}\delta}\left(	\int_{[0,T]^\ell} (D^\ell F)(y){\color{blue}\delta^{k-\ell}}(u({\color{blue}\bullet},y,{\color{red}\bullet}))\d{y} \right)
	\end{equation}
by using the fact that any higher-order divergence operator can be obtained via iteration (cf.\ e.g.\ \cite[formula (2.6.15)]{NouPec12}) and its linearity. Here, and in what follows, the colors correspond to the variables in which the respective operators act. Now, as we have for $\ell\in\{0,\ldots, k\}$ and almost every $y\in [0,T]^\ell$ and $x\in [0,T]$ that $D^\ell F(y)\delta^{k-\ell} (u(\bullet, y,x))\in \D^{1,2}$ with 
	\[ 
		\| D^\ell F(y)\delta^{k-\ell} (u(\bullet, y,x)) \|_{\D^{1,2}}^2\lesssim \supnorm{F}{\ell+1} \|u(\bullet, y,x)\|_{L^2([0,T]^{k-\ell})}^2, 
	\]
we obtain that 
	\[
		\int_{[0,T]^\ell} \int_0^T \| D^\ell F(x)\delta^{k-\ell}(u(\bullet,y,x)) \|_{\D^{1,2}} \d{x} \d{y} \lesssim \supnorm{F}{\ell+1} \|u\|_{L^2([0,T]^{k+1})}<\infty
	\]
so that we can use the Fubini-type theorem for the divergence operator $\delta$ from \cite[Exercise 3.2.7]{Nua06} in \eqref{eq:expansion_1} to obtain the expression
	\begin{equation}
	\label{eq:expansion_2}
		\delta^{k+1}(Fu) =  \sum_{\ell = 0}^k (-1)^\ell {k\choose \ell}  \int_{[0,T]^\ell} {\color{red}\delta} ((D^\ell F)(y){\color{blue}\delta^{k-\ell}}(u({\color{blue}\bullet},y,{\color{red}\bullet})))\d{y}.
	\end{equation} 
Now, as we have, for every $\ell\in \{0,\ldots,k\}$ and almost every $y\in[0,T]^\ell$ that $D^\ell F(y)\in \D^{1,2}$ and ${\color{blue}\delta^{k-\ell}}(u({\color{blue}\bullet},y,{\color{red} \bullet})) \in \Dom {\color{red}\delta}$ with 
	\begin{align*}
		& \E |D^\ell F(y)|^2 \| {\color{blue}\delta^{k-\ell}} (u({\color{blue}\bullet},y,{\color{red} \bullet})) \|_{L^2({\color{red} [0,T] })}^2 \lesssim \supnorm{F}{\ell}^2 \| u({\color{blue}\bullet},y,{\color{red} \bullet})\|^2_{L^2({\color{blue} [0,T]^{k-\ell}}\times {\color{red} [0,T]})} <\infty, \\
		& \E|D^\ell F(y)|^2 |{\color{red} \delta}({\color{blue}\delta^{k-\ell}} (u({\color{blue}\bullet},y,{\color{red} \bullet}))|^2 \lesssim \supnorm{F}{\ell}^2 \|u({\color{blue}\bullet},y,{\color{red}\bullet} )\|_{L^2({\color{blue} [0,T]^{k-\ell}}\times {\color{red} [0,T]})}^2<\infty,\\
		& \E |\langle D^{\ell+1}F(y,{\color{red} \bullet}),{\color{blue}\delta^{k-\ell}}(u ({\color{blue}\bullet},y,{\color{red}\bullet}))\rangle_{L^2({\color{red} [0,T]})}|^2 \lesssim \supnorm{F}{\ell+1}^2 \|u({\color{blue}\bullet},y,{\color{red}\bullet} )\|_{L^2({\color{blue} [0,T]^{k-\ell}}\times {\color{red} [0,T]})}^2<\infty,
	\end{align*}
we can use \cite[Proposition 1.3.3]{Nua06} to factor out the random variable $(D^\ell F)(y)$ from the divergence in \eqref{eq:expansion_2} to obtain 
	\begin{align*}
		\delta^{k+1}(Fu)  & =  \sum_{\ell = 0}^k (-1)^\ell {k\choose \ell}  \int_{[0,T]^\ell} (D^\ell F)(y)\delta^{k+1-\ell}(u(\bullet,y))\d{y} \\
			& \qquad -  \sum_{\ell = 0}^k (-1)^\ell {k\choose \ell} \int_{[0,T]^{\ell+1}} (D^{\ell+1} F)(x) \delta^{k-\ell} (u(\bullet,x))\d{x} \\
			& = \sum_{\ell=0}^{k+1} (-1)^\ell {k+1\choose \ell} \int_{[0,T]^\ell} (D^\ell F)(x) \delta^{k+1-\ell}(u(\bullet,x))\d{x}\\
			& \qquad + \sum_{\ell =1}^k (-1)^{\ell} \left[ {k\choose \ell}-{k+1\choose \ell}\right] \int_{[0,T]^\ell} (D^\ell F)(x)\delta^{k+1-\ell}(u(\bullet,x))\d{x} \\
			& \qquad - \sum_{\ell = 0}^{k-1} (-1)^{\ell} { k \choose \ell} \int_{[0,T]^{\ell+1}} (D^{\ell+1}F)(x)\delta^{k-\ell}(u(\bullet,x))\d{x}
	\end{align*}
and the claim follows by writing $\ell-1$ instead of $\ell$ in the last sum and realizing that 
	\[ 
		{k\choose \ell} -{k+1\choose \ell} + {k\choose \ell-1} = 0, \quad \ell\in \{1,\ldots, k\}.
	\]
\end{proof}

\subsection{Auxiliary lemmas on $p$-variations}

In what follows we give several auxiliary lemmas on $p$-variations that will be needed throughout the proof of \autoref{thm:main}. Recall that we denote 
	\[ V_{n,T}^p(X):= \sum_{i=0}^{2^n-1} |X_{t_{i+1}^n} -X_{t_i^n}|^p\]
for $p\in (0,\infty)$, $n\in\N$, the dyadic partition $\{t_i^n:=Ti2^{-n}\}_{i=0}^{2^n}$ of the interval $[0,T]$, and a stochastic process $X=(X_t)_{t\in [0,T]}$. We begin with a lemma that is easily proved via the Jensen inequality.

\begin{lemma}
\label{lem:triangle_for_var}
Let $n\in\N$ and $p\in [1,\infty)$ and let $X=(X_t)_{t\in [0,T]}$ and $Y=(Y_t)_{t\in [0,T]}$ be two stochastic processes. Then there is the following estimate:
	\[ 
		V_{n,T}^p(X+Y) \leq 2^{p-1}(V_{n,T}^p(X)+V_{n,T}^p(Y)), \quad a.s.
	\]
\end{lemma}

There is also the following lemma that allows to estimate the $L^1$ distance of $p$-variations of two processes. The lemma is proved as \cite[Lemma 4.2]{GueNua05} with $1/H$ therein being replaced by $p$.

\begin{lemma}
\label{lem: est var}
Let $n\in\N$ and $p\in (1,\infty)$, and let $X=(X_t)_{t\in [0,T]}$ and $Y=(Y_t)_{t\in [0,T]}$ be two stochastic processes such that both $\E V_{n,T}^p(X)$ and $\E V_{n,T}^p(Y)$ are finite. Then 
	\[
		\E\left| V_{n,T}^p(X) - V_{n,T}^p(Y)\right|
			\leq p \left(\E V_{n,T}^p(X-Y)\right)^\frac{1}{p} \left( (\E V_{n,T}^p(X))^{1-\frac{1}{p}} + (\E V_{n,T}^p(Y))^{1-\frac{1}{p}}\right).
	\]
\end{lemma}

An application of \autoref{lem: est var} to divergence-type integrals is given now. The lemma generalizes both \cite[Lemma B.2]{CouKri25}, where only deterministic integrands and $p=2$ are considered, and \cite[Lemma 4.3]{GueNua05}, which applies only to an FBM and $p=1/H$.

\begin{lemma}\label{lem: est var int}
    Let $n\in\N$, $H\in (1/2,1)$, let $p\in [1/H,\infty)$, $g,h \in \L^{K,p}:=L^p([0,T];\D^{K,p})$, and let $X=(X_t)_{t\in [0,T]}$ and $Y=(Y_t)_{t\in [0,T]}$ be two stochastic processes defined for $t\in [0,T]$ by 
    \begin{align*}
        X_t = \int_0^t g_s\delta Z^{H,K}_s \quad\mbox{and}\quad  Y_t = \int_0^t h_s\delta Z^{H,K}_s.
    \end{align*}
    Then there is the estimate
    \begin{equation*}
          \E|V_{n,T}^{p}(X) -  V_{n,T}^{p}(Y)| \lesssim \delta_n^{pH-1}\|g-h\|_{\L^{K,p}}(\|g\|^{p-1}_{\L^{K,p}}+\|h\|^{p-1}_{\L^{K,p}}).
    \end{equation*}
\end{lemma}

\begin{proof}
By \autoref{lem:boundedness_of_int_appendix}, it follows that
	\[
		\E V_{n,T}^p(X) = \sum_{i=0}^{2^n-1} \E\left|\int_{t_i^n}^{t_{i+1}^n}g_s\delta{Z^{H,K}_s}  \right|^p \lesssim \sum_{i=0}^{2^n-1}\delta_n^{pH-1} \int_{t_i^n}^{t_{i+1}^n}\|g_s\|_{\D^{K,p}}^p\d{s} = \delta_n^{pH-1}\|g\|_{\L^{K,p}}^p.
	\]
The desired inequality is then obtained by appealing to \autoref{lem: est var}.
\end{proof}

\subsection{Additional auxiliary lemmas for the proof of \autoref{thm:main}}

In what follows, we give the technical \autoref{lem:convergence of integral g} that is essential to the proof of the main \autoref{thm:main}. The following lemma is needed for its proof.

\begin{lemma}
    \label{lem:aux integral ineq}
For arbitrary $T\in (0,\infty)$, $a \in (0,1)$, $\vartheta \geq 0$, and $\varepsilon \in (\vartheta, a + \vartheta)$, we have 
\[
\int_{0}^{T}  \left(\frac{1}{y+x}\right)^{a + \vartheta}  \left(\frac{1}{y}\right)^{1-a}\d{y} \lesssim \left(\frac{1}{x}\right)^\varepsilon, \quad x \in (0,T).
\]
\end{lemma}

\begin{proof}
Note first that there is the convergence
\[
\lim_{x \to 0+}\frac{\int_{0}^{T}  \left(\frac{1}{y+x}\right)^{a + \vartheta}  \left(\frac{1}{y}\right)^{1-a}\d{y}}{\left(\frac{1}{x}\right)^\varepsilon} = \lim_{x \to 0+}\int_{0}^{T}  \left(\frac{x^{\frac{\varepsilon}{a + \vartheta}}}{y+x}\right)^{a + \vartheta}  \left(\frac{1}{y}\right)^{1-a}\d{y} = 0,
\]
where the last equality holds because the integrand converges pointwise to zero as $x\to 0+$ if we find an integrable dominating function. To this end, denote $\delta:= \frac{\varepsilon}{a+ \vartheta} \in (0,1)$ and observe that
\begin{align*}
	\left(\frac{x^\delta}{y+x}\right)^{a + \vartheta}  \left(\frac{1}{y}\right)^{1-a} 
	& \leq \left(\max_{x>0}\left\{\frac{x^\delta}{y+x} \right\}\right)^{a + \vartheta}  \left(\frac{1}{y}\right)^{1-a} \\
	&  =  \left(\delta^\delta (1-\delta)^{(1-\delta)}\frac{1}{y^{1-\delta}}\right)^{a + \vartheta}  \left(\frac{1}{y}\right)^{1-a}\\
	& \lesssim \left(\frac{1}{y}\right)^{1-a + (1-\delta)(a+\vartheta)} \\
	& = \left(\frac{1}{y}\right)^{1+ \vartheta - \varepsilon}
\end{align*}
holds. The last function is integrable on $(0,T)$ as $\varepsilon > \vartheta$. Now, because the function
\[
x \mapsto \frac{\int_{0}^{T}  \left(\frac{1}{y+x}\right)^{a + \vartheta}  \left(\frac{1}{y}\right)^{1-a}\d{y}}{\left(\frac{1}{x}\right)^\varepsilon} 
\]
is continuous on $(0,T]$, the claim of the lemma follows.
\end{proof}

\begin{lemma}
\label{lem:convergence of integral g}
Let $H \in (1/2,1)$, $k\geq 2$ be an integer, and let $\ell\in\{1,2,\ldots, k-1\}$. Then
\begin{equation*}
		\int_{[0,T]^\ell} \left(\int_{0}^{T} |g^\ell (x,u)|^\frac{1}{{H\left(1-\frac{\ell}{k}\right)+ \frac{\ell}{k}}} \d{u}\right)^{H\left(1-\frac{\ell}{k}\right)+\frac{\ell}{k}}\d{x}<\infty,
	\end{equation*}
where 
	\[ g^\ell (x_1,\ldots,x_\ell,u) := c_\ell \prod_{i=1}^\ell \left(\frac{u}{x_i}\right)^{\frac{1}{2}-\frac{1-H}{k}} (u-x_i)_+^{-\left(\frac{1}{2}+\frac{1-H}{k}\right)}
	\]
for $x_1,\ldots, x_\ell\in (0,T]$ and $u\in [0,T]$, in which 
	\[
		 c_\ell := \frac{c_{H,k}}{c_{H\left(1-\frac{\ell}{k}\right)+\frac{\ell}{k},k-\ell}}
	\]
with constant $c_{H,k}$ (and $c_{H\left(1-\frac{\ell}{k}\right) + \frac{\ell}{k},k-\ell}$) being defined by \eqref{eq:cHk}.
\end{lemma}

\begin{proof}
To simplify notation, set 
\[
	a:=\frac{1}{2}-\frac{1-H}{k} \in (0,1/2) 
	\quad \text{ and }\quad 
	\gamma := H\left(1-\frac{\ell}{k}\right)+\frac{\ell}{k} \in (1/2 , 1),
\]
so that we can write
	\[ g^\ell (x_1,\ldots,x_\ell,u) = c_\ell \prod_{i=1}^\ell \left(\frac{u}{x_i}\right)^{a} (u-x_i)_+^{-\left(1-a\right)}
	\]
and
 \begin{align*}
		\mathcal{I}& :=\int_{[0,T]^\ell} \left(\int_{0}^{T} |g^\ell (x,u)|^\frac{1}{{H\left(1-\frac{\ell}{k}\right)+ \frac{\ell}{k}}} \d{u}\right)^{H\left(1-\frac{\ell}{k}\right)+\frac{\ell}{k}}\d{x} \\
        &= \int_{[0,T]^\ell} \left(\int_{0}^{T} \left|c_\ell \prod_{i=1}^\ell \left(\frac{u}{x_i}\right)^{a} (u-x_i)_+^{-\left(1-a\right)}\right|^\frac{1}{\gamma} \d{u}\right)^{\gamma}\d{x}.
	\end{align*}
Denote $\Delta_\ell := \{(x_1,\ldots,x_\ell) \in [0,T]^\ell : x_1 \geq x_2 \geq \ldots \geq x_{\ell}\}$. We have
\begin{align*}
	\mathcal{I} &= \ell! \int_{\Delta_\ell} \left(\int_{x_1}^{T} \prod_{i=1}^\ell \left(\frac{u}{x_i}\right)^{a/\gamma} (u-x_i)^{-\left(1-a\right)/\gamma} \d{u}\right)^{\gamma}\d{x}\\
	&\overset{(1)}{\lesssim} \int_{\Delta_\ell} \prod_{i=1}^\ell \left(\frac{T}{x_i}\right)^{a} \left(\int_{0}^{T-x_1}  \prod_{i=1}^\ell  (v+x_1-x_i)^{-\left(1-a\right)/\gamma} \d{v}\right)^{\gamma}\d{x}\\
	&\overset{}{\lesssim} \int_{\Delta_\ell} \prod_{i=1}^\ell \left(\frac{1}{x_i}\right)^{a} \prod_{i=2}^\ell (x_1-x_i)^{-\left(1-a\right)} \left(\int_{0}^{T}  v^{-\left(1-a\right)/\gamma} \d{v}\right)^{\gamma}\d{x}\\
	&\overset{(2)}{\lesssim} \int_{\Delta_\ell} \prod_{i=1}^\ell \left(\frac{1}{x_i}\right)^{a} \prod_{i=2}^\ell (x_1-x_i)^{-\left(1-a\right)}\d{x}\\
	&\overset{}{=} \int_0^T \int_{x_\ell}^T\ldots \int_{x_3}^T \int_{x_2}^T \prod_{i=1}^\ell \left(\frac{1}{x_i}\right)^{a}  \prod_{i=2}^\ell (x_1-x_i)^{-\left(1-a\right)}\d{x_1}\d{x_2}\ldots \d{x_{\ell-1}}\d{x_\ell}\\
	&\overset{(3)}{=} \int_0^T \int_{x_\ell}^T\ldots \int_{x_3}^T \prod_{i=2}^\ell \left(\frac{1}{x_i}\right)^{a} \int_{0}^{T-x_2}  \left(\frac{1}{y+x_2}\right)^{a}  \prod_{i=2}^\ell (y+ x_2-x_i)^{-\left(1-a\right)}\d{y}\d{x_2}\ldots \d{x_{\ell-1}}\d{x_\ell}\\
	&\overset{}{\lesssim} \int_0^T \int_{x_\ell}^T\ldots \int_{x_3}^T \prod_{i=2}^\ell \left(\frac{1}{x_i}\right)^{a} \prod_{i=3}^\ell (x_2-x_i)^{-\left(1-a\right)} \int_{0}^{T}  \left(\frac{1}{y+x_2}\right)^{a}  \left(\frac{1}{y}\right)^{1-a}\d{y}\d{x_2}\ldots \d{x_{\ell-1}}\d{x_\ell}.
\end{align*}
where in 
	\begin{enumerate}
		\item[(1)] we performed the substitution $v=u-x_1$,
		\item[(2)] we used the fact that $\frac{1-a}{\gamma} = \frac{k+2-2H}{2H(k-\ell)+2\ell}<\frac{k+1}{k+\ell}\leq 1,$
		\item[(3)] we substituted $y=x_1 - x_2$.
	\end{enumerate}
At this point we apply \autoref{lem:aux integral ineq} with $\vartheta =0$ and with $\varepsilon \in (0,a)$ such that $a + (\ell -1) \varepsilon < 1$. By this lemma, we obtain the estimate
\[
	\int_{0}^{T}  \left(\frac{1}{y+x_2}\right)^{a}  \left(\frac{1}{y}\right)^{1-a}\d{y} \lesssim \left(\frac{1}{x_2}\right)^\varepsilon, \quad  x_2 \in (0,T).
\]
This results in
\begin{align*}
	\mathcal{I}&\lesssim \int_0^T \int_{x_\ell}^T\ldots \int_{x_3}^T \prod_{i=2}^\ell \left(\frac{1}{x_i}\right)^{a} \prod_{i=3}^\ell (x_2-x_i)^{-\left(1-a\right)}  \left(\frac{1}{x_2}\right)^\varepsilon \d{x_2}\ldots \d{x_{\ell-1}}\d{x_\ell}\\
&\overset{}{=}  \int_0^T \int_{x_\ell}^T\ldots\int_{x_4}^T\prod_{i=3}^\ell\left(\frac{1}{x_i}\right)^{a} \int_{x_3}^T  \left(\frac{1}{x_2}\right)^{a+\varepsilon} \prod_{i=3}^\ell (x_2-x_i)^{-\left(1-a\right)} \d{x_2} \d{x_3}\ldots \d{x_{\ell-1}}\d{x_\ell}\\
&\overset{}{\lesssim} \int_0^T \int_{x_\ell}^T\ldots\int_{x_4}^T\prod_{i=3}^\ell\left(\frac{1}{x_i}\right)^{a} \prod_{i=4}^\ell (x_3-x_i)^{-\left(1-a\right)} \int_{0}^T  \left(\frac{1}{y+x_3}\right)^{a+\varepsilon}  y^{-\left(1-a\right)} \d{y} \d{x_3}\ldots \d{x_{\ell-1}}\d{x_\ell}.
\end{align*}
Now we can apply \autoref{lem:aux integral ineq} again to obtain
\[
\int_{0}^T  \left(\frac{1}{y+x_3}\right)^{a+\varepsilon}  \left(\frac{1}{y}\right)^{1-a} \d{y}  \lesssim \left(\frac{1}{x_3}\right)^{2 \varepsilon}, \quad x_3 \in (0,T).
\]
Using the above estimate leads to the inequality
\begin{align*}
\mathcal{I}&\lesssim \int_0^T \int_{x_\ell}^T\ldots\int_{x_4}^T\prod_{i=3}^\ell\left(\frac{1}{x_i}\right)^{a} \prod_{i=4}^\ell (x_3-x_i)^{-\left(1-a\right)} \left(\frac{1}{x_3}\right)^{2 \varepsilon} \d{x_3}\ldots \d{x_{\ell-1}}\d{x_\ell}.\\
\end{align*}
By repeating this procedure, we end up with
\begin{align*}
\mathcal{I}&\lesssim \int_0^T \left(\frac{1}{x_\ell}\right)^{a}\int_{x_\ell}^T\left(\frac{1}{x_{\ell-1}}\right)^{a + (\ell-2) \varepsilon}  (x_{\ell-1}-x_\ell)^{-\left(1-a\right)}  \d{x_{\ell-1}}\d{x_\ell}\\
&\overset{}{\lesssim} \int_0^T \left(\frac{1}{x_\ell}\right)^{a}\int_{0}^T\left(\frac{1}{y + x_{\ell}}\right)^{a + (\ell-2) \varepsilon}  \left(\frac{1}{y}\right)^{1-a} \d{x_{\ell-1}}\d{x_\ell}\\
&\overset{}{\lesssim} \int_0^T \left(\frac{1}{x_\ell}\right)^{a + (\ell -1) \varepsilon}\d{x_\ell} < \infty,
\end{align*}
where the convergence of the last integral follows from the choice of $\varepsilon$.
\end{proof}

\end{document}